\documentclass[11pt]{amsart}
\usepackage{amsfonts}

\usepackage{graphicx,calc,amscd}

\newtheorem{theorem}{Theorem}[section]
\newtheorem{corollary}[theorem]{Corollary}
\newtheorem{proposition}[theorem]{Proposition}
\newtheorem{lemma}[theorem]{Lemma}

\newcommand{\Aut}{\operatorname{Aut}}

\newcommand{\Epi}{\operatorname{Epi}}
\newcommand{\G}{\mathcal{G}}

\newcommand{\N}{\mathcal{N}}

\newcommand{\rk}{\func{rk}\,}

\newcommand{\U}{\mathcal{U}}
\newcommand{\V}{\mathbb{V}}

\newcommand{\Z}{\mbox{$\mathbb{Z}$}}
\newcommand{\func}[1]{\text{#1}}

\begin{document}
\title{Orientation-reversing free actions on handlebodies}
\author{Antonio F. Costa}
\address{Departamento de Matematicas Fundamentales\\
Facultad de Ciencias\\
Universidad Nacional de Educaci\'{o}n a Distancia\\
Madrid 28040\\
SPAIN}
\email{acosta@mat.uned.es}
\author{Darryl McCullough}
\address{Department of Mathematics\\
University of Oklahoma\\
Norman, Oklahoma 73019\\
USA}
\email{dmccullough@math.ou.edu}
\urladdr{www.math.ou.edu/$_{\widetilde{\phantom{n}}}$dmccullough/}
\thanks{The first author was supported in part by BFM 2002-04801, the second author
was supported by NSF grant DMS-0102463.}
\date{\today}

\begin{abstract}
We examine free orientation-reversing group actions on orientable
handlebodies, and free actions on nonorientable handlebodies. A
classification theorem is obtained, giving the equivalence classes and weak
equivalence classes of free actions in terms of algebraic invariants that
involve Nielsen equivalence. This is applied to describe the sets of free
actions in various cases, including a complete classification for many (and
conjecturally all) cases above the minimum genus. For abelian groups, the
free actions are classified for all genera.
\end{abstract}

\maketitle

The orientation-preserving free actions of a finite group $G$ on
3-dimension\-al orientable handlebodies have a close connection with a
long-studied concept from group theory, namely \emph{Nielsen equivalence} of
generating sets. The basic result is that the orientation-preserving free
actions of $G$ on the handlebody of genus $g$, up to equivalence, correspond
to the Nielsen equivalence classes of $n$-element generating sets of $G$,
where $n=1+(g-1)/|G|$. This has been known for a long time; it is implicit
in work of J. Kalliongis and A. Miller in the 1980's, as a direct
consequence of theorem~1.3 in their paper \cite{K-Mi1} (for free actions,
the graph of groups will have trivial vertex and edge groups, and the
equivalence of graphs of groups defined there is readily seen to be the same
as Nielsen equivalence on generating sets of $G$). As far as we know, the
first explicit statement detailing the correspondence appears in \cite{MW},
which also contains various applications and calculations using it.

In this paper, we extend the theory from \cite{MW} to free actions that
contain orientation-reversing elements, and to free actions on nonorientable
handlebodies. The orbits of a certain group action on the collection $%
\mathcal{G}_n$ of $n$-element generating sets are the Nielsen equivalence
classes, and this action extends to an action on a set $\mathcal{G}_n\times 
\mathbb{V}_n$, in such a way that the orbits correspond to the equivalence
classes of all free $G$-actions on handlebodies of genus $1+(n-1)|G|$. This
correspondence is our main result, theorem~4.1, which is proven after
presentation of preliminary material on Nielsen equivalence in section~\ref
{sec:Nielsen}, and on ``uniform homeomorphisms'' in section~\ref{sec:uniform}%
. From theorem~4.1, more specific results are derived in section~\ref
{sec:orientable} for orientation-reversing free actions on orientable
handlebodies, and in section~\ref{sec:nonorientable} for free actions on
nonorientable handlebodies. These are illustrated by several calculations
for specific groups, and in section~\ref{sec:abelian} we use the results to
classify all free actions of abelian groups on handlebodies, extending the
classification of orientation-preserving actions given in~\cite{MW}.

We should mention that nonfree actions on handlebodies have been examined in
considerable depth. For nonfree actions, the natural structure on the
quotient object is that of an orbifold, rather than just a handlebody, and
the resulting analysis is much more complicated. A general theory of actions
was given in~\cite{M-M-Z} and the aforementioned~\cite{K-Mi1}, and the
actions on very low genera were extensively studied in~\cite{K-Mi2}. Actions
with the genus small relative to the order of the group were investigated
in~ \cite{M-Z}, and the special case of orientation-reversing involutions is
treated in~\cite{K-M}. The first focus on free actions seems to be~\cite
{Przytycki}, which examines free actions of the cyclic group.

The authors gratefully acknowledge the support of the U. S. National Science
Foundation, the Ministerio de Educaci\'{o}n y Ciencia of Spain, and the
Mathematical Research and Conference Center in Bedlewo of the Polish Academy
of Sciences.

\section{Definitions and notation}

\label{sec:definitions}

In this paper, $G$ will always denote a finite group. A \emph{$G$-action} on
a space $X$ is an injective homomorphism $\Phi\colon G\to \text{Homeo}(X)$.
Two actions $\Phi_1,\Phi_2\colon G\to\text{Homeo}(X)$ are said to be \emph{%
equivalent} if they are conjugate as representations, that is, if there is a
homeomorphism $h\colon X\to X$ such that $h\Phi_1(g)h^{-1}=\Phi_2(g)$ for
each $g\in G$. They are \emph{weakly equivalent} if their images are
conjugate, that is, if there is a homeomorphism $h\colon X\to X$ so that $%
h\Phi_1(G)h^{-1}=\Phi_2(G)$. Equivalently, there is some automorphism $%
\alpha $ of $G$ so that $h\Phi_1(g)h^{-1} = \Phi_2(\alpha(g))$ for all $g$.
In words, equivalent actions are the same after a change of coordinates on
the space, while weakly equivalent actions are the same after a change of
coordinates on the space and a change of the group by automorphism. If $X$
is homeomorphic to $Y$, then the sets of equivalence (or weak equivalence)
classes of actions on $X$ and on $Y$ can be put into correspondence using
any homeomorphism from $X$ to~$Y$.

>From now on, the term \emph{action} will mean a \emph{free} action of a
finite group on a $3$-dimensional handlebody $V_{g}$ of genus $g\geq 1$
(only the trivial group can act freely on the handlebody of genus~$0$, the $%
3 $-ball). One may work in either of the categories of piecewise-linear or
smooth actions. We assume that one of these two categories has been chosen,
and that all maps, isotopies, etc.~lie in that category.

We call an action \emph{orientation-preserving} if $V_{g}$ is orientable and
each element of $G$ acts preserving orientation. We call it \emph{%
orientation-reversing} if $V_{g}$ is orientable and some element of $G$ acts
reversing orientation.

\section{Nielsen equivalence}

\label{sec:Nielsen}

It will be convenient to define Nielsen equivalence in terms of group
actions on sets. We write $C_k$ for the cyclic group of order $k\geq 2$,
including the infinite cyclic group $C_\infty$. Let $\mathbb{U}\cong
C_2*C_2*C_2*C_\infty$ be given by the presentation 
\begin{equation*}
\mathbb{U} = \langle t, u, v, w \;\vert \; t^2=u^2=v^2=1\rangle.
\end{equation*}
For any group $G$ and any positive integer $n\geq 2$, an action of $\mathbb{U%
}$ on the $n$-fold direct product $G^n$ is defined by 
\begin{gather*}
t(g_1,g_2,\ldots,g_n)=(g_1^{-1},g_2,\ldots,g_n) \\
u(g_1,g_2,g_3,\ldots,g_n)=(g_1^{-1},g_1g_2,\ldots,g_n) \\
v(g_1,g_2,g_3,\ldots,g_n)=(g_2,g_1,g_3,\ldots,g_n) \\
w(g_1,g_2,\ldots,g_n)=(g_n,g_1,g_2,\ldots,g_{n-1})\ .
\end{gather*}
\noindent The orbits of this $\mathbb{U}$-action on $G^n$ are called \emph{%
Nielsen equivalence classes.}

Note that if the elements of two Nielsen equivalent $n$-tuples are regarded
as subsets of $G$, then they generate the same subgroup of~$G$. In
particular, if the entries of one of them generate $G$, the same is true for
the other.

Conjugates of $t$ by $w$ allow one to replace any $g_i$ by its inverse.
Conjugates of $v$ by $w$ allow one to interchange any $g_i$ with any $%
g_{i+1} $, and hence to effect any permutation of the coordinates. Simple
combinations of these with $u$ allow one to replace any $g_i$ by $%
g_ig_j^{\pm1}$ or $g_j^{\pm1}g_i$ for some $j\neq i$, keeping all other
coordinates fixed. On the other hand, each of the four generators results
from some sequence of these basic Nielsen ``moves''. Thus Nielsen
equivalence is often described as the equivalence relation generated by
these basic moves.

By letting $\text{Aut}(G)$ act on the left of $G^n$ coordinatewise, we can
extend the $\mathbb{U}$-action to a $\mathbb{U}\times \text{Aut}(G)$-action.
This adds the additional basic Nielsen move 
\begin{equation*}
\alpha(g_1,\ldots,g_n)=(\alpha(g_1),\ldots,\alpha(g_n))
\end{equation*}
for any $\alpha\in \text{Aut}(G)$. The orbits of this $\mathbb{U}\times 
\text{Aut}(G)$-action are called \emph{weak Nielsen equivalence classes}.

We will now see that the Nielsen equivalence classes can be given in terms
of an action of $\text{Aut}(F_n)\times \text{Aut}(G)$ on $G^n$. Defining
this action requires the selection of a basis $x_1,\ldots\,$, $x_n$ of $F_n$%
. Such a selection gives an identification of $G^n$ with the set $\text{Hom}%
(F_n,G)$ of group homomorphisms from $F_n$ to $G$, by regarding $%
(g_1,\ldots, g_n)$ as the homomorphism $\gamma(g_1,\ldots,g_n)\colon F_n\to
G $ that sends $x_i$ to $g_i$. The action of $\text{Aut}(F_n)\times \text{Aut%
}(G)$ on $G^n$ is then defined simply by $(\phi,\alpha)\cdot
\gamma=\alpha\circ \gamma \circ \phi^{-1}$.


We regard the restriction of this action to $\text{Aut}(F_n)\times\{1\}$ as
an $\text{Aut}(F_n)$-action. The next lemma shows that the action of $%
\mathbb{U}\times \text{Aut}(G)$ on $G^n$ always factors through the action
of $\text{Aut}(F_n)\times \text{Aut}(G)$ on $G^n$.

\begin{lemma}
The orbits of the $\text{Aut}(F_{n})$-action on $G^{n}$ (respectively, the $%
\text{Aut}(F_{n})\times \text{Aut}(G)$-action on $G^{n}$) are exactly the
Nielsen equivalence classes (respectively, the weak Nielsen equivalence
classes). In fact, there is a surjective homomorphism $A_{n}\colon \mathbb{U}%
\rightarrow \text{Aut}(F_{n})$ such that the action of an element $(u,\alpha
)\in \mathbb{U}\times \text{Aut}(G)$ equals the action of $(A_{n}(u),\alpha )
$. Changing the choice of basis for $F_{n}$ changes $A_{n}$ by an inner
automorphism of $\text{Aut}(F_{n})$. \label{lem:antihomomorphism}
\end{lemma}

\begin{proof}
Define $T\in \Aut(F_n)$ by $T(x_1)=x_1^{-1}$ and $T_i(x_j)=x_j$ for $j>1$,
and similarly define $U$, $V$, and $W$ corresponding to $u$, $v$, and
$w$. It is straightforward to check that
$(t,\alpha)(g_1,\ldots,g_n)=(T,\alpha)(g_1,\ldots, g_n)$, and similarly for
the other three generators, so the action of $\mathbb{U}$ on $G^n$ factors
through the image of the ``capitalization'' function $A_n\colon
\mathbb{U}\to \Aut(F_n)$.  Using well-known generating sets for
$\Aut(F_n)$, such as that of Nielsen's presentation \cite{Nielsen} or the
Fouxe-Rabinovitch presentation listed in \cite{MM}, one checks that $A_n$
is surjective. The basis change remark is a straightforward check.
\end{proof}

We will now extend Nielsen equivalence in $G^n$ to a relation on a larger
set that will capture some orientation information when we apply it to study
actions on handlebodies.

Write $\mathbb{V}_{n}$ for the direct sum $\displaystyle\oplus_{i=1}^{n}C_{2}
$, where $C_{2}=\{-1,1\}$. Using the selected basis $x_{1},\ldots \,$, $x_{n}
$ of $F_{n}$, identify $\mathbb{V}_{n}$ with $\text{Hom}(F_{n},C_{2})$ by
identifying an element $(v_{1},\ldots ,v_{n})$ of $\mathbb{V}_{n}$ with the
homomorphism $\omega (v_{1},\ldots ,v_{n})$ that sends $x_{i}$ to $v_{i}$.
We define an $\text{Aut}(F_{n})\times \text{Aut}(G)$-action on $G^{n}\times 
\mathbb{V}_{n}$ by putting 
\begin{equation*}
(\phi ,\alpha )\cdot (\gamma ,\omega )=(\alpha \circ \gamma \circ \phi
^{-1},\omega \circ \phi ^{-1})\ .
\end{equation*}
Restricted to the subset $G^{n}\times \{(1,\ldots ,1)\}$, this can be
identified with the $\text{Aut}(F_{n})\times \text{Aut}(G)$-action on $G^{n}$%
.

We use $\mathcal{G}_n$ to denote the set of generating $n$-vectors of $G$,
that is, $n$-tuples $(g_1,\ldots,g_n)$ of elements of $G$ such that $%
\{g_1,\ldots,g_n\}$ generates~$G$. These correspond to the surjective
elements of $\text{Hom}(F_n,G)$, so $\mathcal{G}_n\times \mathbb{V}_n$ is
invariant under the $\text{Aut}(F_n)\times \text{Aut}(G)$-action.

\section{Uniform homeomorphisms}

\label{sec:uniform}

We will use an idea which has appeared several times in the literature \cite
{CM}, \cite{Montreal}, \cite{MM} (the most relevant of these references is 
\cite{MM}, since it also concerns handlebodies). The quotient of a free
action on a genus $g$ handlebody is a handlebody $V_{n}$ of genus~$%
n=1+(g-1)/|G|$ (see section~\ref{sec:classification}). This handlebody is
regarded as one component of a disjoint union of a family of handlebodies
indexed by $\mathbb{V}_{n}$, where the handlebody $N(v_{1},\ldots ,v_{n})$
corresponding to a vector $(v_{1},\ldots ,v_{n})$ has the property that
traveling around the $i^{th}$ handle reverses the local orientation exactly
when $v_{i}=-1$. An $n$-tuple $(g_{1},\ldots ,g_{n})$ of elements that
generate $G$ determines a $G$-action on a handlebody with quotient $%
N(v_{1},\ldots ,v_{n})$ in the following way: $G$ acts by covering
transformations on the covering space of $N(v_{1},\ldots ,v_{n})$
corresponding to the kernel of the homomorphism $\pi _{1}(N(v_{1},\ldots
,v_{n}))\rightarrow G$ that sends the generator corresponding to the $i^{th}$
handle to $g_{i}$.

A key property of this family of handlebodies is that any element of $\text{%
Aut}(\pi_1(V_n))$ can be realized, in an appropriate sense, by a ``uniform''
homeomorphism of the family. The action of uniform homeomorphisms on the set
of components of the family corresponds exactly to the $\text{Aut}(F_n)$%
-action on $\mathbb{V}_n$ defined in section~\ref{sec:Nielsen}. Uniform
homeomorphisms overcome the technical problem that an automorphism of $%
\pi_1(V_n)$ need not preserve the orientability of $1$-handles and hence
need not be induced by a self-homeomorphism of $V_n$.

The proof of the main technical result, theorem~\ref{thm:classification},
shows that two pairs $((g_1,\ldots,\allowbreak g_n),(v_1,\ldots, v_n))$ and $%
((g_1^{\prime},\ldots, g_n^{\prime}),(v_1^{\prime},\ldots, v_n^{\prime}))$
in $\mathcal{G}_n\times \mathbb{V}_n$ lie in the same $\text{Aut}(F_n)$%
-orbit exactly when there is a homeomorphism between $N(v_1,\ldots, v_n)$
and $N(v_1^{\prime},\ldots, v_n^{\prime})$ that lifts to an equivalence
between the actions which have them as quotients and are determined by $%
(g_1,\ldots,g_n)$ and $(g_1^{\prime},\ldots, g_n^{\prime})$.

Here is the construction from \cite{MM}. Fixing a positive integer $n$, let $%
R_n$ be a $1$-point union of $n$ circles. Write $F_n$ for the free group $%
\pi_1(R_n)$. Let $x_1,\ldots\,$, $x_n$ be the standard set of generators of $%
F_n$, where $x_i$ is represented by a loop that travels once around the $%
i^{th}$ circle.

To set notation, let $\Sigma$ be a $3$-ball, and in $\partial \Sigma$ select 
$2n$ disjoint imbedded $2$-disks $D_1,E_1, D_2,E_2,\ldots,D_n,E_n$. Fix
orientation-preserving imbeddings $J_i\colon D^2\to D_i$ and $K_i\colon
D^2\to E_i$. Let $r\colon D^2\to D^2$ send $(x,y)$ to $(x,-y)$. For $%
v=(v_1,\ldots,v_n)\in \mathbb{V}_n$, construct a handlebody $N(v)$ as
follows. For each $i$, let $H_i$ be a copy of $D^2\times I$ and identify $%
(x,y,0)$ with $J_i(x,y)$ and $(x,y,1)$ with $K_ir^{(1+v_i)/2}(x,y)$. The
resulting $1$-handle $H_i$ is orientation-preserving or
orientation-reversing according as $v_i$ is $1$ or~$-1$.

Regard $N(v)$ as a thickening of $R_n$, in which the join point is the
center $*$ of $\Sigma$ and the loop of $R_n$ that represents $x_i$ goes once
over $H_i$ from $D_i$ to $E_i$ and does not meet any other $H_j$. Traveling
around this $i^{th}$ loop preserves the local orientation at $*$ if and only 
$v_i=1$. Thus $N(1,\ldots,1)$ is orientable, while all other $N(v)$ are
nonorientable and are homeomorphic to $N(-1,\ldots,-1)$. We denote the
disjoint union of the $N(v)$ by~$\mathcal{N}$.

We will now define a homeomorphism of $\mathcal{N}$ called a uniform slide
homeomorphism. For each $N(v)$, write $N^{\prime }(v)$ for the closure of $%
N(v)-H_{1}$. Choose a loop $\alpha $ in $\partial N^{\prime }(v)$, based at
the origin in $E_{1}$, that travels through $\partial \Sigma $ to $\partial
E_{2}$, once over $H_{2}$ to $\partial D_{2}$, and returns in $\partial
\Sigma $ to the origin of $E_{1}$. There is an isotopy $J_{t}$ of $N^{\prime
}(v)$ such that

\begin{enumerate}
\item  $J_{0}$ is the identity of $N^{\prime }(v)$,

\item  each $J_{t}$ the identity outside a regular neighborhood of $%
E_{1}\cup \alpha $,

\item  during $J_{t}$, $E_{1}$ moves once around $\alpha $, traveling over $%
H_{2}$ from $E_{2}$ to $D_{2}$, and

\item  the restriction of $J_{1}$ to $E_{1}$ is the identity or $r$,
according as $J_{1}$ preserves or reverses the local orientation on $E_{1}$.
\end{enumerate}

\noindent A homeomorphism of $\mathcal{N}$ is defined by sending $N(v)$ to $%
N(w)$ using $J_{1}$ on $N^{\prime }(v)$ and the identity on $H_{1}$. Here, $%
(w_{1},\ldots ,w_{n})=(v_{1}v_{2},v_{2},\ldots ,v_{n})$, since the $r$ in
item (4) will be needed exactly when $v_{2}=-1$. There are many choices of
sliding loop $\alpha $, nonisotopic in $\partial N^{\prime }(v)$, so the
homeomorphism of $\mathcal{N}$ is by no means uniquely defined up to isotopy.

With respect to the identifications $\pi _{1}(R_{n})=\pi _{1}(N(v))$ given
by the inclusions of $R_{n}$ into $N(v)$ and $N(w)$, the homeomorphism from $%
N(v)$ to $N(w)$ induces the automorphism $\rho $ of $F_{n}$ that sends $%
x_{1} $ to $x_{1}x_{2}$ and fixes all other $x_{j}$. Note that $%
(w_{1},w_{2},\ldots ,w_{n})=\rho \cdot (v_{1},v_{2},\ldots ,v_{n})$, for the
action of $\rho $ on $(v_{1},v_{2},\ldots ,v_{n})$ defined in section~\ref
{sec:Nielsen}.

This particular basic slide homeomorphism is called sliding the right end
(that is, $E_{1}$) of $H_{1}$ over $H_{2}$. Similarly, one can uniformly
slide the right or left end of any $H_{i}$ over any other $H_{j}$, either
from $E_j$ to $D_j$ or from $D_j$ to $E_j$, obtaining homeomorphisms whose
effect on components of $\mathcal{N}$ agrees with the action of their
induced automorphisms on $\mathbb{V}_{n}$. These are called \emph{uniform
slide homeomorphisms} of~$\mathcal{N}$.

A \emph{uniform interchange} of $H_i$ and $H_j$ is defined using an isotopy $%
J_t$ that interchanges both $D_i$ and $D_j$, and $E_i$ and $E_j$. It sends $%
N(\ldots,v_i,\ldots,v_j,\ldots)$ to $N(\ldots,v_j,\ldots,v_i,\ldots)$, and
induces the automorphism of $F_n$ that interchanges $x_i$ and $x_j$. Using a 
$J_t$ that interchanges $D_i$ and $E_i$ defines a \emph{uniform spin} of the 
$i^{th}$ handle. This preserves each component of $\mathcal{N}$, and induces
the automorphism that sends $x_i$ to $x_i^{-1}$.

There are two other kinds of basic uniform homeomorphisms, both of which
preserve each $N(v)$ and induce the identity automorphism on $F_n$. Choose a
reflection of $\Sigma$ that preserves $*$ and restricts to $r$ on each $D_i$
and $E_i$. Define a homeomorphism of $N(v)$ by taking $r\times 1_I$ on each $%
H_i$ and the chosen reflection on $\Sigma$. The resulting uniform
homeomorphism of $\mathcal{N}$ is denoted by $R$. Finally, any Dehn twist
about a properly imbedded $2$-disk in $\mathcal{N}$ is a basic uniform
homeomorphism.

In all cases, the action of the basic uniform homeomorphism on the
components of $\mathcal{N}$ agrees with the action on $\mathbb{V}_n$ of the
automorphism it induces on $F_n$ with respect to the identifications $%
F_n=\pi_1(R_n)=\pi_1(N(v))$.

A \emph{uniform homeomorphism} of $\mathcal{N}$ is a homeomorphism (freely)
isotopic to a composition of the basic uniform homeomorphisms we have
defined here. The inverse of a basic uniform homeomorphism is a basic
uniform homeomorphism, so the inverse of any uniform homeomorphism is
uniform.

By abuse of notation, we write $*$ for the union of the basepoints of the
components of $\mathcal{N}$, and by $\mathcal{M}(\mathcal{N},*)$ the group
of isotopy classes of homeomorphisms of $\mathcal{N}$ that preserve this
subset. The uniform homeomorphisms that preserve $*$ form a subgroup $%
\mathcal{U}(\mathcal{N},*)$ of $\mathcal{M}(\mathcal{N},*)$, called the 
\emph{uniform mapping class group.} We mention that although we have given
infinitely many generators, it can be shown that $\mathcal{U}(\mathcal{N},*)$
is finitely generated. This is proven in \cite{MM}.

For $v\in \mathbb{V}_n$, let $\text{St}(N(v),*)\subseteq \mathcal{U}(%
\mathcal{N},*)$ be the stabilizer of the component $N(v)$ under the action
of $\mathcal{U}(\mathcal{N},*)$ on the components of $\mathcal{N}$. We have
the following result from~\cite{MM}:

\begin{theorem}
The restriction $\text{St}(N(v),\ast )\rightarrow \mathcal{M}(N(v),\ast )$
is surjective. Any homeomorphism $N(v)\rightarrow N(w)$ is isotopic to the
restriction of a uniform homeomorphism. \label{thm:stabilizer}
\end{theorem}

\begin{proof} The first statement is basically theorem~7.2.3 
from~\cite{MM}, proven there for compression bodies, which include
handlebodies as a special case. The restriction in \cite{MM} to mapping
classes of local degree~$1$ at $*$ is not needed since we have included the
reflection $R$ among our uniform homeomorphisms.

For the second statement, note first that the uniform homeomorphisms act
transitively on the set of nonorientable components of $\mathcal{N}$, so
given $g\colon N(v)\to N(w)$, there is a uniform homeomorphism $u_1$ that
carries $N(w)$ to $N(v)$. (To see this, suppose that $N(w)$ and $N(v)$ are
nonorientable and choose some $w_i=-1$. Slide the other handles of $N(w)$
over the $i^{th}$ handle as necessary to make $w_j=v_j$ for $j\neq i$. If
all these $w_j$ are now $1$, then $w_i=-1=v_i$ since $N(w)$ and $N(v)$ are
nonorientable. If not, there is some other $w_j=-1$, and a slide of the
$i^{th}$ handle over the $j^{th}$ can be used if needed to change $w_i$ to
equal $v_i$.)  By the first sentence of the theorem, the composition
$u_1\circ g$ is isotopic to the restriction of a uniform homeomorphism
$u_2$ that stabilizes $N(v)$, so on $N(v)$, $g$ is isotopic to
$u_1^{-1}\circ u_2$.
\end{proof}

\section{The algebraic classification of actions}

\label{sec:classification}

Suppose that $\Phi\colon G\to \text{Homeo}(V)$ is a free action on a
handlebody $V$, possibly nonorientable. Its quotient $N$ is also a
handlebody. To see this, recall that any torsionfree finite extension of a
finitely generated free group is free (by~\cite{K-P-S} any finitely
generated virtually free group is the fundamental group of a graph of groups
with finite vertex groups, and if the group is torsionfree, the vertex
groups must be trivial), so $\pi_1(V/G)$ is free. Since $V$ is irreducible,
so is $V/G$, and theorem~5.2 of~\cite{Hempel} shows that $V/G$ is a
handlebody.

>From covering space theory, the action $\Phi$ determines an extension 
\begin{equation*}
1\longrightarrow \pi_1(V)\longrightarrow \pi_1(N)\overset{\pi}{%
\longrightarrow} G\longrightarrow 1
\end{equation*}
where $\pi(x)$ is defined by taking a representative loop for $x$, lifting
it to a path starting at the basepoint of $V$, and letting $\pi(x) $ be the
covering transformation that sends the basepoint of $V$ to the endpoint of
the path. Writing $n$ for the genus of $N$, the Euler characteristic shows
that $1+|G|\,(n-1)$ is the genus of~$V$. The genus of $N $ can be any $n$
greater than or equal to $\mu(G)$, the minimum number of elements in a
generating set of $G$. In particular, the genera of handlebodies on which $G$
acts freely preserving orientation are exactly $1+\vert G \vert (n-1)$ where 
$n\geq \mu(G)$. The \emph{minimal genus} is $1+\vert G\vert (\mu(G)-1)$.

Choose any $N(v_{1},\ldots ,v_{n})$ that is homeomorphic to $N$, and choose
a homeomorphism $k\colon N\rightarrow N(v_{1},\ldots ,v_{n})$. Let $W$ be
the covering of $N(v_{1},\ldots ,v_{n})$ determined by the subgroup $%
k_{\#}(\pi _{1}(V))$ (where $k_{\#}$ is the isomorphism induced by $k$ on
the fundamental groups). This subgroup is well-defined up to conjugacy, so $W
$ depends only on the choice of $k$. Let $K\colon V\rightarrow W$ be a lift
of $k$. It identifies $G$ with the group of covering transformations of $W$.
Each basis element $x_{i}$ of $F_{n}=\pi _{1}(N(v_{1},\ldots ,v_{n}))$
(where the $x_{i}$ are as in section~\ref{sec:uniform}) determines a
covering transformation $g_{i}\in G$. We associate to $\Phi $ the pair $%
((g_{1},\ldots ,g_{n}),(v_{1},\ldots ,v_{n}))$, which we will abbreviate as $%
(g,v)$. Since the $x_{i}$ generate $F_{n}$, the $g_{i}$ generate $G$, so $%
(g,v)$ is an element of $\mathcal{G}_{n}\times \mathbb{V}_{n}$.

\begin{theorem}
Sending $\Phi $ to the orbit of the element $(g,v)$ defines a bijection from
the equivalence classes (respectively, weak equivalence classes) of free $G$%
-actions on handlebodies of genus $1+|G|\,(n-1)$ to the set of $\text{Aut}%
(F_{n})$-orbits (respectively, $\text{Aut}(F_{n})\times \text{Aut}(G)$%
-orbits) in $\mathcal{G}_{n}\times \mathbb{V}_{n}$. \label%
{thm:classification}
\end{theorem}

\begin{proof}
Changing the choice of basepoint in $W$ or the lift of $k$ changes
$((g_1,\ldots,g_n),v)$ to
$((hg_1h^{-1},\ldots,hg_nh^{-1}),v)$ for some $h\in G$.
Choose an element $\widetilde{h}\in F_n$ with 
$\gamma(g_1,\ldots,g_n)(\widetilde{h})=h$,
and let $\mu(\widetilde{h})\in \Aut(F_n)$ be the automorphism that
conjugates by $\widetilde{h}^{-1}$. Then
$((hg_1h^{-1},\ldots,hg_nh^{-1}),v)=
(\mu(\widetilde{h}),1)\cdot((g_1,\ldots,g_n),v)$, so these
elements lie in the same $\Aut(F_n)$-orbit.

Suppose a different $N(v')$ and $k'\colon N\to N(v')$ are used to associate
a pair $(g',v')=((g_1',\ldots,g_n'),(v_1',\ldots,v_n'))$ to~$\Phi$. By
theorem~\ref{thm:stabilizer}, $k'\circ k^{-1}\colon N(v)\to N(v')$ is the
restriction of a uniform homeomorphism~$u$. We claim that $(u_\#,1)\in
\Aut(F_n)$ carries $(g,v)$ to $(g',v')$. Since the action of $\U(\N,*)$ on
the components of $\N$ induces the action of $\Aut(F_n)$ on $\V_n$, it
suffices to show that $\gamma(g_1,\ldots,g_n)\circ
u_\#^{-1}=\gamma(g_1',\ldots,g_n')$, that is, that
$g_i'=\gamma(g_1,\ldots,g_n)(u_\#^{-1}(x_i))$.

Let $(W,w)$ and $(W',w')$ be the covering spaces of $N(v_1,\ldots,v_n)$ and
$N(v_1',\ldots,v_n')$ respectively, such that lifting $x_i$ to $W$ and $W'$
produces $g_i$ and $g_i'$ respectively. Let $\widetilde{u}\colon (W,w)\to
(W',w')$ be the lift of $u|_{N(v)}$. Now, $g_i'$ is the covering
transformation that carries $w'$ to the endpoint of the lift of $x_i$
starting at $w'$. Consider $(u|_{N(v)})^{-1}(x_i)$. Its lift to $W$
starting at $w$ is carried by $\widetilde{u}$ to the lift of $x_i$ in $W'$
starting at $w'$.  That is, the covering transformation of $W$
corresponding to $g_i'$ under $\widetilde{u}$ is determined by
$u_\#^{-1}(x_i)$, so is $\gamma(g_1,\ldots,g_n)(u_\#^{-1}(x_i))$. This
verifies the claim.

Equivalent actions produce equivalent associated elements. For if $\Phi$ is
equivalent to another $G$-action $\Phi'$ on $V'$, with quotient $N'$, then
there is a homeomorphism $j\colon N'\to N$ that lifts to an equivariant
homeomorphism from $V'$ to $V$. Since we may use $k\circ j$ as the
homeomorphism from $N'$ to $N(v)$ to define the element associated to
$\Phi'$, the associated pairs are in the same $\Aut(F_n)$-orbit.

Conversely, suppose that the pairs $(g,v)$ and $(g',v')$ associated to the
actions $\Phi$ and $\Phi'$ are in the same $\Aut(F_n)$-orbit. Let $\phi\in
\Aut(F_n)$ carry one to the other. By theorem~\ref{thm:stabilizer}, there
is a uniform homeomorphism $u\in \U(\N)$ inducing $\phi$, which must
carry $N(v_1,\ldots,v_n)$ to $N(v_1',\ldots,v_n')$. The condition that
$\gamma\circ \phi^{-1}=\gamma'$ ensures that $u$ lifts to a $G$-equivariant
homeomorphism from $(W,w)$ to $(W',w')$, so the actions on these covering
spaces are equivalent.  Since the actions on $W$ and $W'$ are respectively
equivalent to the original actions on $V$ and $V'$, the original actions
were equivalent.

Finally, being able to apply an automorphism of $G$ at any point in the
process changes equivalence to weak equivalence, and enlarges the choices
of $(g,v)$ to the $\Aut(F_n)\times \Aut(G)$-orbit.
\end{proof}

Under the action of $\text{Aut}(F_n)$ (or $\text{Aut}(F_n)\times \text{Aut}%
(G)$) on $\mathbb{V}_n$ the element $(1,\ldots,1)$ is fixed, so the subset $%
\mathcal{G}_n\times\{(1,\ldots,1)\}$ is a union of orbits, in fact the
orbits correspond exactly to the Nielsen equivalence classes (or weak
equivalence classes) of elements of $\mathcal{G}_n$. This recovers the
algebraic classification of orientation-preserving actions in 
\cite[Theorem~2.3]{MW}.

\section{Actions on orientable handlebodies}

\label{sec:orientable}

>From theorem~\ref{thm:classification}, an explicit representative of the
equivalence class of $G$-actions corresponding to the $\text{Aut}(F_n)$%
-orbit of the element $(g,v)$ of $\mathcal{G}_n\times \mathbb{V}_n$ is the
covering space $W$ of $N(v)$ whose fundamental group is the kernel of $%
\gamma=\gamma(g_1,\ldots,g_n)\colon F_n\to G$. Since $v_i$ tells the
orientability of $x_i$ in $N(v)$, a covering space is orientable if and only
if it corresponds to a subgroup in the kernel of $\omega=\omega(v_1,%
\ldots,v_n)$. Therefore there is a simple criterion for $W$ to be orientable:

\begin{proposition}
Let $W$ be the covering space of $N(v)$ corresponding to the kernel of $%
\gamma $. Then $W$ is orientable if and only if there is $\overline{\omega }%
\in \text{Hom}(G,C_{2})$ such that $\omega \colon F_{n}\rightarrow C_{2}$
factors as $\overline{\omega }\circ \gamma \colon F_{n}\rightarrow
G\rightarrow C_{2}$. Equivalently, sending $g_{i}$ to $v_{i}$ defines a
homomorphism from $G$ to $C_{2}$. \label{prop:orientability}
\end{proposition}

Applying theorem~\ref{thm:classification}, we obtain:

\begin{corollary}
Under the correspondence of theorem~\ref{thm:classification}, the
equivalence classes (respectively, weak equivalence classes) of free $G$%
-actions on orientable handlebodies of genus $1+|G|\,(n-1)$ correspond to
the set of $\text{Aut}(F_{n})$-orbits (respectively, $\text{Aut}%
(F_{n})\times \text{Aut}(G)$-orbits) in $\mathcal{G}_{n}\times \mathbb{V}_{n}
$ for which sending $g_{i}$ to $v_{i}$ (on one, hence on any representative)
determines a homomorphism $\overline{\omega }$ from $G$ to~$C_{2}$. \label%
{coro:orientation-reversing}
\end{corollary}

It will be useful to make explicit the induced action of $\text{Aut}%
(F_n)\times \text{Aut}(G)$ on these $\overline{\omega}$. In the statement of
proposition~\ref{prop:epiaction}, we call $\overline{\omega}$ the element of 
$\text{Hom}(G,C_2)$ \emph{associated to $(g,v)$.}

\begin{proposition}
If $\overline{\omega }\in \text{Hom}(G,C_{2})$ is associated to $(g,v)\in 
\mathcal{G}_{n}\times \mathbb{V}_{n}$ and $(\phi ,\alpha )\in $Aut$%
(F_{n})\times $Aut$(G)$, then $\overline{\omega }\circ \alpha ^{-1}$ is the
element of $\text{Hom}(G,C_{2})$ associated to $(\phi ,\alpha )\cdot (g,v)$.

\label{prop:epiaction}
\end{proposition}

\begin{proof}
Regarding $(g,v)$ as $(\gamma, \omega)$ we have $(\phi,\alpha)\cdot
(\gamma,\omega) = (\alpha\circ \gamma\circ \phi^{-1},\omega\circ
\phi^{-1})$. Since $\omega\circ \phi^{-1}= (\overline{\omega}\circ
\alpha^{-1}) \circ (\alpha\circ \gamma\circ \phi^{-1})$, its associated
element is $\overline{\omega}\circ \alpha^{-1}$.
\end{proof}

The classification up to equivalence of free actions on orientable
handlebodies is no more difficult than the classification of generating $n$%
-vectors of $G$ up to Nielsen equivalence. For $n\geq \mu(G)$ let $\mathcal{E%
}_n$ denote the set of Nielsen equivalence classes of generating $n$-vectors
of $G$. We write $\text{Epi}(G,C_2)$ for the set of surjective homomorphisms
from $G$ to $C_2$, that is, all elements of $\text{Hom}(G,C_2)$ except the
trivial homomorphism~$0$.

\begin{theorem}
For $n\geq \mu (G)$, the set of equivalence classes of free $G$-actions on
the orientable handlebody of genus $1+|G|(n-1)$ corresponds bijectively to $%
\mathcal{E}_{n}\times \text{Hom}(G,C_{2})$, with the orientation-preserving
actions corresponding to $\mathcal{E}_{n}\times \{0\}$ and the
orientation-reversing actions corresponding to $\mathcal{E}_{n}\times \text{%
Epi}(G,C_{2})$. \label{thm:orientable}
\end{theorem}

\begin{proof}
By theorem~\ref{thm:classification}, every action is equivalent to the
action of $G$ by covering transformations on a covering space $W$ of some
$N(v)$, and the equivalence classes of actions correspond to the
$\Aut(F_n)$-orbits of $\G_n\times \V_n$. Restricting to the
$\G_n$-coordinate defines a function $\G_n\times \V_n \to \G_n$ which is
$\Aut(F_n)$-equivariant, so there is an induced function on the sets of
$\Aut(F_n)$-orbits.

Fix an $\Aut(F_n)$-orbit of $\G_n$ and a generating $n$-vector
$(h_1,\ldots\,h_n)$ that represents it. Each $\Aut(F_n)$-orbit of
$\G_n\times \V_n$ that restricts to this element contains a representative
of the form $((h_1,\ldots,h_n),(v_1,\ldots,v_n))$. The element
$((h_1,\ldots,h_n),(1,\ldots,1))$ is not equivalent to any other such
element, and represents the unique element that corresponds to an
orientation-preserving action. By
corollary~\ref{coro:orientation-reversing},
$((h_1,\ldots,h_n),(v_1,\ldots,v_n))$ corresponds to an
orientation-reversing action if and only if sending $h_i$ to $v_i$ defines
a surjective homomorphism from $G$ to $C_2$. By
proposition~\ref{prop:epiaction}, this homomorphism is an invariant of the
equivalence class. On the other hand, each element $\omega$ of
$\Epi(G,C_2)$ determines a choice of $v$ for which $\omega=\omega(v)$, so
the equivalence classes of orientation-reversing actions that restrict to
the orbit of $(h_1,\ldots,h_n)$ in $\G_n$ correspond to $\Epi(G,C_2)$.
\end{proof}

For classification of orientation-reversing actions up to weak equivalence,
there is an added difficulty. An $\text{Aut}(F_n)\times \text{Aut}(G)$-orbit
of elements of $\mathcal{G}_n$ is a union of a collection of $\text{Aut}%
(F_n) $-orbits, say $\{C_1,\ldots,C_r\}$. It produces one weak equivalence
class of orientation-preserving actions, but for orientation-reversing
actions, one must determine the $\text{Aut}(G)$-orbits of $%
\{C_1,\ldots,C_r\}\times \text{Epi}(G,C_2)$. This seems to be a subtle
problem, in general.

It often happens, however, that $\mathcal{G}_n$ consists of only one $\text{%
Aut}(F_n)$-orbit, in which case the action of $\text{Aut}(G)$ on $%
\{C_1\}\times \text{Epi}(G,C_2)$ can be identified with the action on $\text{%
Epi}(G,C_2)$. Thus in this case, the classification of actions on the
orientable handlebody $V_g $ is easy:

\begin{theorem}
Suppose that all elements of $\mathcal{G}_{n}$ are Nielsen equivalent, and
put $g={1+|G|(n-1)}$. Then

\begin{enumerate}
\item  There is only one equivalence class of orientation-preserving free $G$%
-actions on $V_{g}$.

\item  The set of weak equivalence classes of orientation-reversing free
actions of $G$ on $V_{g}$ corresponds bijectively to the set of $\text{Aut}%
(G)$-orbits of $\text{Epi}(G,C_{2})$.
\end{enumerate}

\label{thm:weakorientable}
\end{theorem}

Conjecturally, all generating $n$-vectors are equivalent whenever $G$ is
finite and $n>\mu(G)$ (see the discussion in \cite{MW}). So the previous
theorem might give a complete classification of all actions on orientable
handlebodies above the minimal genus. The conjecture is known for many
classes of groups, such as solvable groups \cite{Dunwoody1}, $\text{PSL}%
(2,p) $ ($p$ prime) \cite{Gilman}, $\text{PSL}(2,3^p)$ ($p$ prime) \cite{MW}%
, $\text{PSL}(2,2^m)$ \cite{Evans1}, and the Suzuki groups $\text{Sz}%
(2^{2m-1})$ \cite{Evans1}.

A nice example is the quaternion group $Q$ of order $8$. One can check that
for any $n\geq 2=\mu(Q)$, any two generating $n$-vectors of $Q$ are Nielsen
equivalent. So for any $k\geq 1$, there is one equivalence class of
orientation-preserving free $Q$-action on $V_{1+8k}$, and there are three
equivalence classes of orientation-reversing free $Q$-actions, corresponding
to the nonzero elements of $\text{Hom}(Q,C_2)=H^1(Q;\text{\mbox{$\mathbb{Z}$}%
}/2)=\text{\mbox{$\mathbb{Z}$}}/2\oplus \text{\mbox{$\mathbb{Z}$}}/2$. Under
the $\text{Aut}(Q)$-action on $\text{Epi}(Q,C_2)$, all three elements lie in
the same orbit, so there is only one weak equivalence class of
orientation-reversing free $Q$-action on $V_{1+8k}$.

Let us finish this section with another example. For $r\geq 3$ let $D_{r}$
be the dihedral group of $2r$ elements and presentation $\left\langle
s_{1},s_{2}:s_{1}^{2}=s_{2}^{2}=(s_{1}s_{2})^{r}=1\right\rangle$.

Suppose first that $n>2=\mu (D_{r})$. Since $D_{r}$ is solvable, there is
only one $\text{Aut}(F_{n})$-orbit in $\mathcal{G}_{n}$, and hence there is
only one equivalence class of orientation-preserving actions. If $r$ is
even, there are three classes of orientation-reversing actions, represented
by the elements 
\begin{gather*}
\{((s_{1},s_{2},1,...,1),(-1,-1,1,...,1)),
((s_{1},s_{2},1,...,1),(1,-1,1,...,1)), \\
((s_{1},s_{2},1,...,1),(-1,1,1,...,1))\}
\end{gather*}
of $\mathcal{G}_n\times \text{Epi}(D_r,C_2)$. If $r$ is odd, the second two
do not define homomorphisms from $D_{r}$ to $C_2$, and there is only one
equivalence class. When $r$ is even, there are two weak equivalence classes
of orientation-reversing actions, represented by: 
\begin{equation*}
\{((s_{1},s_{2},1,...,1),(-1,-1,1,...,1)),
((s_{1},s_{2},1,...,1),(1,-1,1,...,1))\}.
\end{equation*}
Suppose now that $n=2$. A set of representatives of the $\text{Aut}(F_{2})$%
-orbits in $\mathcal{G}_{2}$ is $\{(s_{1},(s_{1}s_{2})^{m}):1\leq
m<r/2,(m,r)=1\}$ (see theorem~4.5 of \cite{MW}). There are $\varphi (r)/2$
classes of orientation-preserving actions (where $\varphi$ is the Euler
function) forming one weak equivalence class. If $r$ is is odd there are $%
\varphi (r)/2$ classes of orientation-reversing actions. If $r$ is even
there are $3\varphi (r)/2$ classes of orientation-reversing actions forming $%
\varphi (r)$ weak equivalence classes.

\section{Actions on nonorientable handlebodies}

\label{sec:nonorientable}

There is a simple algebraic criterion for $G$ to act freely on the
nonorientable handlebody $N_m$ of genus $m$. Recall that $H^1(G;\text{%
\mbox{$\mathbb{Z}$}}/2)$ can be identified with $\text{Hom}(G,C_2)$.

\begin{proposition}
$G$ acts freely on $N_{m}$ if and only if $m=1+|G|(n-1)$ where $n\geq \mu (G)
$ and $n>\text{rk}\,H^{1}(G;\mbox{$\mathbb{Z}$}/2)$. \label%
{prop:nonorientablecriterion}
\end{proposition}

\begin{proof}
If $n<\mu(G)$ then $\G_n$ is empty
and $G$ does not act freely on any handlebody of genus $n$,
so we assume that $n\geq \mu(G)$. According to
corollary~\ref{coro:orientation-reversing}, an element
$((h_1,\ldots,h_n),(v_1,\ldots,v_n))\in \G_n\times \V_n$ represents an
orbit corresponding to an action on a nonorientable handlebody if and only
if sending $h_i$ to $v_i$ does not define a homomorphism from $G$ to
$C_2$. So $N_m$ has no free action exactly when all of the $2^n$ choices
for $v$ define homomorphisms. Since $\rk H^1(G;\Z/2)\leq \mu(G)\leq n$, 
the latter is equivalent to $\rk H^1(G;\Z/2) = n$.
\end{proof}

For the quaternion group $Q$ considered in section~\ref{sec:orientable}, we
have $2=\mu(Q)=\text{rk}\, H^1(Q;\text{\mbox{$\mathbb{Z}$}})$, so $Q$ acts
freely on $V_9$, but not on~$N_9$.

We may combine proposition~\ref{prop:nonorientablecriterion} with theorem~%
\ref{thm:orientable} to determine the genera on which $G$ can act:

\begin{corollary}
Let $A=\{1+|G|(n-1)\;|\;n\geq \mu (G)\}$. Then

\begin{enumerate}
\item  $G$ acts freely preserving orientation on $V_{m}$ if and only if $%
m\in A$.

\item  $G$ acts freely reversing orientation on $V_{m}$ if and only if $m\in
A$ and $\text{rk}\,H^{1}(G;\mbox{$\mathbb{Z}$}/2)>0$.

\item  $G$ acts freely on $N_{m}$ if and only if $m\in A$ and either $%
m>1+|G|(\mu (G)-1)$ or $\text{rk}\,H^{1}(G;\mbox{$\mathbb{Z}$}/2)<\mu (G)$.
\end{enumerate}

\label{coro:genera}
\end{corollary}

There is a version of theorem~\ref{thm:weakorientable} for actions on
nonorientable handlebodies.

\begin{theorem}
Suppose that $n>\mu (G)$ and that all generating $n$-vectors of $G$ are
Nielsen equivalent. Put $m={1+|G|(n-1)}$. Then, all free actions of $G$ on $%
N_{m}$ are equivalent. \label{thm:redundant}
\end{theorem}

\begin{proof}
Fix a generating set $h_1,\ldots\,$, $h_{n-1}$ with $n-1$ elements. Since
all generating $n$-vectors are Nielsen equivalent, each $\Aut(F_n)$-orbit
of $\G_n\times \V_n$ has a representative of the form $((h_1,\ldots,
h_{n-1},1),v)$. Fix such an element corresponding to an action on~$N_m$.

Suppose first that $v_n=-1$. For any $i$ with $v_i=1$, the basic Nielsen
move sending $h_i$ to $h_ih_n=h_i$ changes $v_i$ to $v_iv_n=-1$. So
$((h_1,\ldots, h_{n-1},1),v)$ is equivalent to $((h_1,\ldots,
h_{n-1},1),(-1,\ldots,-1))$.  Suppose that $v_n=1$.  By
corollary~\ref{coro:orientation-reversing}, sending each $h_i$ to $v_i$
does not define a homomorphism to $C_2$, so there is some product
$h_{i_1}^{\epsilon_1}\cdots h_{i_k}^{\epsilon_k}=1$, with all
$\epsilon_i=\pm1$, for which $v_{i_1}^{\epsilon_1}\cdots
v_{i_k}^{\epsilon_k}=-1$.  A sequence of $k$ basic Nielsen moves replacing
$h_n$ by $h_nh_{i_j}^{\epsilon_j}$ shows that $((h_1,\ldots,
h_{n-1},1),(v_1,\ldots,v_{n-1},1))$ is equivalent to $((h_1,\ldots,
h_{n-1},1),(v_1,\ldots,v_{n-1},-1))$, which we have seen is equivalent to
$((h_1,\ldots, h_{n-1},1),(-1,\ldots,-1))$. Therefore we have only one
$\Aut(F_n)$-orbit of elements of $\G_n\times \V_n$ that corresponds to an
action on a nonorientable handlebody.
\end{proof}

By way of illustration, we return to our example of actions of $D_{r}$. If $r
$ is odd then $\text{rk}\,H^{1}(D_{r},\mbox{$\mathbb{Z}$}_{2})=1$ and if $r$
is even then $\text{rk}\,H^{1}(D_{r},\mbox{$\mathbb{Z}$}_{2})=2$. We have $%
\mu (D_{r})=2$, and corollary~\ref{coro:genera} shows that $D_{r}$ acts on $%
N_{2r+1}$ if and only if $r$ is odd. When $r$ is odd, there are $\varphi
(m)/2$ equivalence classes of actions, represented by $%
((s_{1},(s_{1}s_{2})^{m}),(-1,-1))$ where $m$ is relatively prime to $r$ and 
$1\leq m<r/2$. These form one weak equivalence class. When $n>2$, there is
one equivalence class of actions on $N_{1+2r(n-1)}$ represented by $%
((s_{1},s_{2},1,\ldots ,1),(-1,-1,-1,\ldots ,-1))$.

As we noted in section~\ref{sec:orientable}, it is conjectured that all
generating $n$-vectors are equivalent whenever $G$ is finite and $n>\mu(G)$,
so theorem~\ref{thm:redundant} might classify all actions on $N_m$ when $%
m>1+\vert G\vert (\mu(G)-1)$. The classification of actions on the
nonorientable handlebody of genus $1+\vert G\vert (\mu(G)-1)$ seems to be an
interesting general problem.

\section{Actions of abelian groups}

\label{sec:abelian}

In this section, we will completely classify free actions of abelian groups
on handlebodies.

Throughout this section, we assume that $G$ is abelian. For now, write $G$
as $C_{d_1}\oplus\cdots \oplus C_{d_n}$ where $d_{i+1}|d_i$ for $1\leq i<n$.
We have $\mu(G)=n$, since clearly $\mu(G)\leq n$, while $G\otimes
C_{d_n}\cong C_{d_n}^n$ requires $n$ generators.

Theorem~4.1 of \cite{MW} tells the equivalence classes of generating $\mu
(G) $-vectors. Fix a generator $s_{i}$ for $C_{d_{i}}$. Each $\text{Aut}%
(F_{n})$-orbit in $\mathcal{G}_{n}\times \mathbb{V}_{n}$ contains exactly
one element of the form $(s_{1},\ldots ,s_{n-1},s_{n}^{m})$ where $m$ is
relatively prime to $d_{n}$ and $1\leq m\leq d_{n}/2$. There is only one
weak equivalence class, since for each such $m$, there is an automorphism of 
$G$ fixing $s_{i}$ for $i<n$ and sending $s_{n}$ to $s_{n}^{m}$.

It will be convenient to rewrite $G$ as $C_{e_1}\oplus \cdots \oplus
C_{e_k}\oplus C_{d_1}\oplus \cdots \oplus C_{d_\ell}$, where the $e_i$ are
even, the $d_j$ are odd, each $e_{i+1}\vert e_i$, each $d_{j+1}\vert d_j$,
and $d_1\vert e_k$. We write $s_i$ for the selected generator of $C_{e_i}$
and $t_j$ for the selected generator of $C_{d_j}$. There is a corresponding
decomposition $\mathbb{V}_n=\mathbb{V}_k\oplus \mathbb{V}_\ell$, in which we
will denote elements by $(v,w)=(v_1,\ldots,v_k,w_1,\ldots,w_\ell)$. Also, we
write $|\{e_1,\ldots,e_k\}|$ for the cardinality of the set $%
\{e_1,\ldots,e_k\}$.

We now analyze the $\text{Aut}(F_n)$- and $\text{Aut}(F_n)\times \text{Aut}%
(G)$-orbits on $\mathcal{G}_n\times \mathbb{V}_n$. Using theorem~4.1 of \cite
{MW} discussed above, every $\text{Aut}(F_n)$-orbit has a representative of
the form $((s_1,\ldots,s_k,t_1,\ldots,t_\ell^m),(v_1,\ldots,v_k,w_1,\ldots,
w_\ell))$, or of the form $((s_1,\ldots,s_k^m),(v_1,\ldots,v_k))$ if $\ell=0$%
. For such a representative, choose a corresponding free action of $G$ on a
handlebody $W$.

Suppose first that $W$ is orientable. Proposition~\ref{prop:orientability}
shows that all $w_j=1$. Each choice of $v$ determines a different
homomorphism $\overline{\omega}\colon G\to C_2$, so all the possible choices
for $v$ (an element of $\mathbb{V}_k$) and $m$ (an integer relatively prime
to $d_\ell$ with $1\leq m\leq d_\ell/2$, or relatively prime to to $e_k$
with $1\leq m\leq e_k/2$ if $\ell=0$) determine inequivalent actions. As in
theorem~\ref{thm:orientable}, the choices with $v=(1,\ldots,1)$ are the
orientation-preserving actions, and all others are orientation-reversing.

Still assuming that $W$ is orientable, we consider weak equivalence. If $%
\alpha$ is the automorphism of $G$ that sends $t_\ell$ to $t_\ell^m$ (or $%
s_k $ to $s_k^m$, when $\ell=0$) then the action of $(1,\alpha)$ sends $%
((s_1,\ldots,s_k,t_1,\ldots,t_\ell),(v_1,\ldots,v_k,1,\ldots,1))$ to $%
((s_1,\ldots,s_k,t_1,\ldots,t_\ell^m),(v_1,\ldots,v_k,1,\ldots, 1))$ (or $%
((s_1,\ldots,s_k),(v_1,\ldots,v_k))$ to $((s_1,\ldots,s_k^m),(v_1,%
\ldots,v_k))$), so for weak equivalence we may eliminate the orbit
representatives with $m\neq 1$. In particular, there is only one weak
equivalence class of orientation-preserving actions. Suppose the action is
orientation-reversing, so that some $v_j=-1$. Choose the largest such $j$.

Suppose that $v_{i}=1$ for some $e_{i}$ for which $e_{j}|e_{i}$. Let $\alpha 
$ be the automorphism of $G$ that sends $s_{i}$ to $s_{i}s_{j}$ and fixes
all other generators, and let $\rho $ be the automorphism of $F_{n}$ that
sends $x_{i}$ to $x_{i}x_{j}$ and fixes all other generators. We have $(\rho
,\alpha )\cdot ((s_{1},\ldots ,s_{k},t_{1},\ldots ,t_{\ell }),(v_{1},\ldots
,v_{i},\ldots ,v_{j},\ldots v_{k},1,\ldots ,1))=((s_{1},\ldots
,s_{k},t_{1},\ldots ,t_{\ell }),(v_{1},\ldots ,v_{i}v_{j},\ldots
,v_{j},\ldots v_{k},1,\ldots ,1))$. Repeating this for all such $i$, we may
make $v_{i}=-1$ whenever $e_{j}|e_{i}$; that is, after possibly reselecting $%
j$ to a larger value with the same value of $e_{j}$, we may assume that $%
v_{i}=-1$ for every $i\leq j$, $v_{i}=1$ for every $i>j$, and that $%
e_{j+1}<e_{j}$ (or $j=k$). Taking only representatives with this property
reduces our collection of representatives of $\text{Aut}(F_{n})\times \text{%
Aut}(G)$-orbits to only $|\{e_{1},\ldots ,e_{k}\}|$ elements. To check that
no two of these can be in the same orbit, we observe that the kernels of the 
$\overline{\omega }$ for these different elements are not isomorphic.
Alternatively we may think in terms of actions: For the action defined by an
element in this form, there is a primitive element in $\pi _{1}(N(v))$ that
determines an orientation-reversing covering transformation of $W$, and
whose $e_{j}$-th power lifts to an orientation-preserving loop, and $e_{j}$
is the smallest integer with this property. For every action weakly
equivalent to this one, $e_{j}$ must be the smallest integer with this
property.

Suppose now that $W$ is nonorientable, and again consider an orbit
representative $((s_1,\ldots,s_k,t_1,\ldots,t_\ell^m),(v_1,\ldots,v_k,w_1,%
\ldots, w_\ell))\in \mathcal{G}_n\times \mathbb{V}_n$. Proposition~\ref
{prop:orientability} shows that some $w_j=-1$. By basic Nielsen moves
replacing an $s_i$ (or a $t_i$) by $s_it_j$ (or $t_it_j$) $d_j$ times, we
may make every $v_i$ and every $w_i$ equal to $-1$ (in case $j=\ell$, use $%
t_\ell^m$ rather than $t_\ell$). Therefore the equivalence classes of
actions correspond to the choices for $m$, and there is only one weak
equivalence class.

We now collect these observations.

\begin{theorem}
Let $G=C_{e_{1}}\oplus \cdots \oplus C_{e_{k}}\oplus C_{d_{1}}\oplus \cdots
\oplus C_{d_{\ell }}$, as above. If $e_{k}=2$, put $N=1$, otherwise put $%
N=\varphi (e_{k})/2$ if $\ell =0$ and $N=\varphi (d_{\ell })/2$ if $\ell >0$%
. Then the free actions on handlebodies of minimal genus $1+|G|(k+\ell -1)$
are as follows.

\begin{enumerate}
\item  For orientation-preserving actions, there are $N$ equivalence
classes, forming one weak equivalence class.

\item  For orientation-reversing actions, there are $(2^{k}-1)N$ equivalence
classes, forming $|\{e_{1},\ldots ,e_{k}\}|$ weak equivalence classes.

\item  If $\ell =0$, then $G$ does not act freely on the nonorientable
handlebody. If $\ell >0$, then there are $N$ equivalence classes, forming
one weak equivalence class.
\end{enumerate}

\label{thm:abelian}
\end{theorem}

For actions above the minimal genus, we have:

\begin{theorem}
For $n>k+\ell $, $G$ acts freely on the orientable and nonorientable
handlebodies of genus $1+|G|(n-1)$, with the following equivalence classes.

\begin{enumerate}
\item  For orientation-preserving actions, there is one equivalence class.

\item  For orientation-reversing actions, there are $2^{k}-1$ equivalence
classes, forming $|\{e_{1},\ldots ,e_{k}\}|$ weak equivalence classes.

\item  For actions on the nonorientable handlebody, there is one equivalence
class.
\end{enumerate}

\label{thm:higherabelian}
\end{theorem}

\begin{proof}
Since $G$ is solvable, \cite{Dunwoody1} shows that all generating
$n$-vectors Nielsen are equivalent to $(s_1,\ldots,
s_k,t_1,\ldots,t_\ell,1,\ldots,1)$. Therefore theorem~\ref{thm:orientable}
gives part~(1) and theorem~\ref{thm:redundant} gives~(3). For~(2), the
proof is then essentially the same as that of theorem~\ref{thm:abelian}; if
one allows some of the $d_i$ to equal $1$, in effect making $k+\ell=n$,
then the proof is almost line-for-line unchanged.
\end{proof}

\end{document}